\documentclass[a4paper]{article}
\usepackage{amsmath,amssymb,amsfonts,amsthm}
\usepackage[ngerman]{babel}
\usepackage[T1]{fontenc}
\usepackage[latin9]{inputenc}
\usepackage{enumitem}
\usepackage{tikz}
\usepackage{geometry}

\def\d{\partial}
\def\In{\subseteq}
\renewcommand\l{\lambda}
\def\cM{\mathcal{M}}
\def\MfG{\cM_f(G)}
\def\N{\mathbb{N}}
\def\cP{\mathcal{P}}
\def\PfG{\cP_f(G)}
\renewcommand\phi{\varphi}
\def\R{\mathbb{R}}
\def\cV{\mathcal{V}}

\def\D{\Delta}
\def\Df{\D_f}
\def\Dff{\Df^*}
\def\w{{\scriptstyle{\mathcal{W}}}}
\def\wf{\w_f}
\def\wff{\wf^*}
\def\ci{\chi'}
\def\cif{\ci_f}
\def\ciff{\ci^{\,*}_f}
\def\Gammaf{\Gamma_f}
\def\Gammaff{\Gammaf^*}

\newcommand{\set}[2]{\{#1\,|\;#2\}}
\def\all{\forall\,}
\def\lflr{\left\lfloor}
\def\rflr{\right\rfloor}

\theoremstyle{plain}
\newtheorem{theorem}{Theorem}[]
\newtheorem{corollary}[]{Corollary}

\newtheorem{conjecture}{Conjecture}
\newtheorem{lemma}[]{Lemma}

\theoremstyle{remark}
\newtheorem{remark}{Remark}
\newtheorem{example}{Example}

\newenvironment{myproof}{\noindent{\em Proof.}}{\hfill{\hbox{\rule{2mm}{2mm}}}\par\medskip}

\begin{document}

\title{\textbf{On the $f$-matching polytope and the fractional $f$-chromatic index\\[11pt]}}

\author{Stefan Glock\\[11pt]}
\date{}

\maketitle
\renewcommand*{\abstractname}{Abstract}
\begin{abstract}
Our motivation is the question of how similar the $f$-colouring problem is to the classic edge-colouring problem, particularly with regard to graph parameters. 
In 2010, Zhang, Yu, and Liu \cite{Zhang} gave a new description of the $f$-matching polytope and derived a formula for the fractional $f$-chromatic index, stating that the fractional $f$-chromatic index equals the maximum of the fractional maximum $f$-degree and the fractional $f$-density. Unfortunately, this formula is incorrect. We present counterexamples for both the description of the $f$-matching polytope and the formula for the fractional $f$-chromatic index. Finally, we prove a short lemma concerning the generalization of Goldberg's conjecture. 
\medskip

%

\end{abstract}

\section{Introduction}

Throughout this paper, the term graph refers to a finite and undirected graph, which may have multiple edges but no loops. The vertex set and the edge set of a graph $G$ are denoted by $V(G)$ and $E(G)$, respectively. 
If $X$ and $Y$ are subsets of $V(G)$, then $E_G(X,Y)$ contains all edges that connect $X$ and $Y$. Let $E_G[X]$ denote the set of all edges with both ends in $X$ and $\d_GX$ the set of all edges with exactly one end in $X$. Thus, the degree of a vertex $v$ in $G$ is $d_G(v)=|\d_G\{v\}|$. If the meaning is clear from the context, we will frequently omit superfluous subscripts and brackets for the sake of readability. For example, henceforth, we will write $\d v$ instead of $\d_G\{v\}$. The expression $H\In G$ means $H$ is a subgraph of $G$, and for $U\In V(G)$, the induced subgraph is denoted by $G[U]$.

A weighted graph is a pair $(G,f)$ consisting of a graph $G$ and a vertex function $f$ of $G$, which assigns a positive integer to every vertex of $G$. For $U\In V(G)$, set $f(U)=\sum_{v\in U}f(v)$, and $f(G)$ should stand for $f(V(G))$. An $f$-matching of the weighted graph $(G,f)$ is an edge set $M\In E(G)$ so that each vertex $v\in V(G)$ satisfies $|M\cap\,\d v|\leq f(v)$. The set of all $f$-matchings of $G$ is denoted by $\MfG$.

An $f$-colouring, introduced by Hakimi and Kariv \cite{Hakimi}, assigns to every edge of $G$ a colour, satisfying that at each vertex $v$ each colour occurs at most $f(v)$ times. More formal, $\phi\colon E(G)\rightarrow C$ is an $f$-colouring of $G$ iff $\phi^{-1}(\alpha)\in\MfG$ for all $\alpha\in C$, where $C$ is an arbitrary set. The $f$-chromatic index, denoted by $\cif(G)$, is the least possible cardinality of such a colour set. As the computation of $\cif$ is NP-complete, one is interested in good bounds. The mere fact that every $f$-colouring induces a partition of $E(G)$ into $f$-matchings gives rise to two easy lower bounds, the maximum $f$-degree and the $f$-density. Let us first define $$\Dff(G)=\max_{v\in V(G)}\frac{d_G(v)}{f(v)}$$ as the fractional maximum $f$-degree and $$\wff(G)=\max_{H\In G,\;|V(H)|\geq 2}\frac{|E(H)|}{\lflr\frac{1}{2}f(H)\rflr}$$ as the fractional $f$-density, where we set $\wff(G)=0$ if $G$ has less than two vertices. Then, the maximum $f$-degree of $G$ is $\Df(G)=\lceil\Dff(G)\rceil$ and the $f$-density of $G$ is defined by $\wf(G)=\lceil\wff(G)\rceil$. Easy observation yields
\begin{equation} \label{lower bound}
\cif(G)\geq\max\{\Df(G),\wf(G)\}\,\mbox{.}
\end{equation}
It is much more complicated to find good upper bounds. In 1988, Nakano, Nishizeki, and Saito \cite{Nakano} proved that any weighted graph satisfies $$\cif(G)\leq\max\left\{\frac{9}{8}\Df(G)+\frac{6}{8},\wf(G)\right\}\,\mbox{,}$$ which encouraged them to transfer Goldberg's conjecture to the $f$-colouring problem.

\begin{conjecture} \label{Nakano}

Any weighted graph satisfies $$\cif(G)\leq\max\{\Df(G)+1,\wf(G)\}\,\mbox{.}$$

\end{conjecture}
If this proves to be true the $f$-chromatic index would be restricted to the values $\Df(G),\Df(G)+1$ and $\wf(G)$. Of course, the computation of $\wf(G)$ seems to be NP-hard as well, however, the value of $\max\{\Df(G)+1,\wf(G)\}$ can be computed efficiently. This observation is closely linked to the fractional $f$-chromatic index, which can be defined in several ways. We want to do it by means of fractional $f$-colourings. 

A fractional $f$-colouring of $G$ is a map $w\colon\MfG\rightarrow[0,1]$ satisfying the following condition:
\begin{equation} \label{condition}
\sum_{M\in\MfG\colon e\in M}w(M)=1\qquad\all e\in E(G)
\end{equation}

For a fractional $f$-colouring $w$ of $G$, we call $$\sum_{M\in\MfG}w(M)$$ the value of $w$. The fractional $f$-chromatic index $\ciff(G)$ is then the minimum value over all fractional $f$-matchings of $G$, which exists, since this is an LP-problem bounded from below. Note that if one replaced the closed interval $[0,1]$ by the set $\{0,1\}$, the obtained minimum value would be nothing else than the $f$-chromatic index $\cif(G)$, where the function $w$ indicates whether a certain $f$-matching is a colour class or not. Thus, the fractional $f$-chromatic index is a lower bound for the $f$-chromatic index.
\begin{remark}
One could replace $[0,1]$ by the nonnegative real numbers and '$=$' by '$\geq$' in (\ref{condition}) in order to define fractional $f$-colourings, as Zhang et al. \cite{Zhang} did. That does not change the obtained minimum value (see \cite{Stiebitz}, Theorem B.1).
\end{remark}

\begin{remark}
While the computation of the $f$-chromatic index is NP-complete, the fractional $f$-chromatic index can be determined efficiently. We refer the reader to \cite{Scheinerman} and \cite{Schrijver} for more profound information on algorithmic details and computational complexity.
\end{remark}

\section{The $f$-matching polytope}

Let $(G,f)$ be an arbitrary weighted graph with at least one edge. Let $\cV(G)$ then denote the real vector space of all functions $x\colon E(G)\rightarrow\R$, which is isomorphic to the standard vector space $\R^{|E(G)|}$. The characteristic function of an edge set $F\In E(G)$ is denoted by $i_F$, where
$$i_F(e)=
\begin{cases}
1 & \mbox{if }e\in F,\\
0 & \mbox{if }e\not\in F.
\end{cases}
$$

The $f$-matching polytope $\PfG$ of $G$ is then defined as the convex hull of the characteristic functions of all $f$-matchings, i.e., $$\PfG=\mathrm{conv}\,(\set{i_M}{M\in\MfG}).$$ If $f(v)=1$ for all vertices of $G$, we write $\cP(G)$ instead of $\PfG$, which stands for the ordinary matching polytope.

The $f$-matching polytope is defined over its extreme points which is admittedly impractical. We are now interested in a description by a system of linear inequalities, which exists in any case. Edmonds \cite{Edmonds} was the first who accomplished that. For a vector $x\in\cV(G)$ and an edge set $F\In E(G)$, define $x(F)=\sum_{e\in F}x(e)$.

\begin{theorem} \label{Edmonds1}

For any graph $G$ with at least one edge, a vector $x\in\cV(G)$ belongs to the matching polytope iff $x$ satisfies the following system of linear inequalities:

    \begin{enumerate}[label=\rm{(\arabic*)}]

        \item   $\all e\in E(G)\colon x(e)\geq0$

        \item   $\all v\in V(G)\colon x(\d v)\leq1$

        \item   $\all U\In V(G)\colon x(E[U])\leq\lflr\frac{1}{2}|U|\rflr$

    \end{enumerate}

\end{theorem}

\begin{remark}
One can easily find that the characteristic function of any matching satisfies (1) -- (3), and therefore every convex combination of these functions satisfies the system. The crucial statement is that the given inequalities suffice to determine the matching polytope. However, not all of them are really necessary, for instance, in (3) it would be enough to include subsets $U$ with odd cardinality.
\end{remark}

Edmonds \cite{Edmonds} also gave a description of the $f$-matching polytope. Note that this statement does not immediately imply Theorem~\ref{Edmonds1} for $f\equiv1$.

\begin{theorem} \label{Edmonds2}

Let $G$ be an arbitrary graph with at least one edge and $f$ a vertex function of $G$. A vector $x\in\cV(G)$ belongs to the $f$-matching polytope iff $x$ satisfies the following system of linear inequalities:

    \begin{enumerate}[label=\rm{(\alph*)}]

        \item   $\all e\in E(G)\colon 0\leq x(e)\leq1$

        \item   $\all v\in V(G)\colon x(\d v)\leq f(v)$

        \item   $\all U\In V(G)\,\all F\In \d U\colon x(E[U]\cup F)\leq\lflr\frac{1}{2}(f(U)+|F|)\rflr$

    \end{enumerate}

\end{theorem}


In order to proof their description of the $f$-matching polytope, Zhang et al. \cite{Zhang} embraced the strategy of a proof for Theorem~\ref{Edmonds1} given by Scheinerman and Ullman \cite{Scheinerman}. For a weighted graph $(G,f)$, let $Q_f(G)$ denote the polyhedron in $\cV(G)$ described by the following inequalities:

\begin{enumerate}[label=\rm{(\roman*)}]

    \item   $\all e\in E(G)\colon x(e)\geq0$

    \item   $\all v\in V(G)\colon x(\d v)\leq f(v)$

    \item   $\all U\In V(G)\colon x(E[U])\leq\lflr\frac{1}{2}f(U)\rflr$

\end{enumerate}

For $f\equiv1$, this is exactly the system given by Theorem~\ref{Edmonds1}. The assertion of Zhang, Yu, and Liu \cite{Zhang} reads: For any weighted graph $(G,f)$, the $f$-matching polytope $\PfG$ is equal to $Q_f(G)$. One may mention that they consider only bounded vertex functions, i.e., $f(v)\leq d(v)$ for all $v\in V(G)$. Although this actually does not make any difference, we will construct a counterexample with bounded vertex function.

\begin{example}

There is a very simple counterexample showing that $\PfG\neq Q_f(G)$ in general. Consider the graph $G$ consisting of two vertices that are connected by two edges and a vertex function $f$ that assigns the value $2$ to both vertices. The function $x\in\cV(G)$ may map one edge on $2$ and the other edge on $0$. Then $x$ satisfies (i) -- (iii), but does obviously not belong to the $f$-matching polytope.

\end{example}

The problem we encountered here is apparently that there is no constraint $x(e)\leq1$ in the definition of $Q_f(G)$. In Theorem~\ref{Edmonds1} this constraint can be omitted, for it follows directly from (2). However, we need it for weighted graphs. Moreover, Zhang et al. actually need it in their proof of \textsc{Claim B} to have reasonable evidence that they `always can find a $\Pi$-extremal $\mathbf{i}_{F'}$ such that $e\not\in F'$'.\footnote{Zhang et al.\cite{Zhang}, p. 3364} Since we are anyhow mainly interested in the inequalities (c) and (iii), respectively, we just add the constraint $\all e\in E(G)\colon x(e)\leq1$ and call the obtained polyhedron once again $Q_f(G)$. Yet the following counterexample will show that neither this can save the statement.

\begin{example}

Let $k\geq3$ be an odd integer. Consider the graph $G$ comprising a cycle of length $k$ and an extra vertex $u'$ that is connected with the vertex $u$ on the cycle by two edges $e_1,e_2$. The vertex function $f$ may assign the value $2$ to $u$ and $u'$, and $1$ to the other vertices. Note that $f(v)\leq d(v)$ for all $v\in V(G)$. Let us define the vector $x\in\cV(G)$ by

$$x(e)=\begin{cases}
1 &\mbox{if }e=e_1,\\
0 &\mbox{if }e=e_2,\\
\frac{1}{2} &\mbox{otherwise}.\end{cases}$$

We now have to verify that $x$ satisfies (i) -- (iii), but is not contained in $\PfG$.

\medskip
\begin{myproof}
That $x$ satisfies (i) results directly from its definition. With $x(\d u)=1+0+\frac{1}{2}+\frac{1}{2}\leq2,x(\d u')=1+0\leq2$, and $x(\d v)=\frac{1}{2}+\frac{1}{2}\leq1$ for all other vertices $v\in V(G)$, we also have verified (ii). In order to check (iii), we distinguish some cases. Let $U$ be an arbitrary subset of $V(G)$. If both $u$ and $u'$ are contained in $U$, then both $e_1$ and $e_2$ belong to $E[U]$ and we therefore have $x(E[U])=1+0+\frac{1}{2}(|E[U]|-2)=\frac{1}{2}|E[U]|$. Otherwise, neither $e_1$ nor $e_2$ belongs to $E[U]$ and this also yields $x(E[U])=\frac{1}{2}|E[U]|$. We show now that $|E[U]|\leq f(U)-1$ holds. 

\smallskip
\noindent{\textbf{Case 1.} $\quad u\not\in U$}
\begin{itemize}
\item $G[U]$ 
does not contain a cycle $\quad\Rightarrow\quad |E[U]|\leq |U|-1$
\item $f(U)\geq |U|$
\end{itemize}
\noindent{\textbf{Case 2.} $\quad u\in U,\,u'\not\in U$}
\begin{itemize}
\item $2|E[U]|=\sum_{v\in U}d_{G[U]}(v)\leq|U|\cdot2\quad\Rightarrow\quad|E[U]|\leq|U|$
\item $f(U)=1\cdot2+(|U|-1)\cdot1=|U|+1$
\end{itemize}
\noindent{\textbf{Case 3.} $\quad u,u'\in U$}
\begin{itemize}
\item $2|E[U]|=\sum_{v\in U}d_{G[U]}(v)\leq1\cdot4+1\cdot2+(|U|-2)\cdot2\quad\Rightarrow\quad|E[U]|\leq|U|+1$
\item $f(U)=2\cdot2+(|U|-2)\cdot1=|U|+2$
\end{itemize}
In each case, the mentioned inequality holds, and we can deduce $$x(E[U])=\frac{1}{2}|E[U]|\leq\frac{1}{2}(f(U)-1)\leq\lflr\frac{1}{2}f(U)\rflr.$$ Thus, $x$ satisfies (iii).

Now, assume that $x$ belongs to the $f$-matching polytope of $G$. Consequently, we have $$x=\sum_{M\in\MfG}\l_Mi_M\quad\mbox{for suitable}\quad\l_M\in\R_{\geq0}\quad\mbox{with}\quad\sum_{M\in\MfG}\l_M=1\,\mbox{.}$$ For an edge $e\in E(G)$, we obtain $$x(e)=\sum_{M\in\MfG}\l_Mi_M(e)=\sum_{M\in\MfG\colon e\in M}\l_M\,\mbox{.}$$ Now, we add up both sides of this equation over the edges of the cycle of $G$. On the left side, we easily obtain $\frac{k}{2}$. Since $x(e_1)=1$, all $\l_M$ with $e_1\not\in M$ need to vanish. And since every $f$-matching that contains $e_1$ can contain at most $\lfloor\frac{k}{2}\rfloor$ edges from the cycle, we infer that each $\l_M\neq0$ appears at most $\lfloor\frac{k}{2}\rfloor$-times on the right side. As the sum of all $\l_M$ equals $1$, the value of the right side is at most $\lfloor\frac{k}{2}\rfloor$. Thus, we have shown $\frac{k}{2}\leq\lfloor\frac{k}{2}\rfloor$, a contradiction to our premise that $k$ is odd. Ergo, $x$ can not belong to the $f$-matching polytope of $G$.
\end{myproof}

\end{example}

\begin{remark}
One may ask, at which point the proof of Zhang et al. is flawed. Without introducing all needed terminology, we just want to mention that the problem appears in {\it Subcase \sl 3.2 }, when they consider an $F F'$-alternating walk $Q$ in $G-e$ starting from $u$. This walk could be extended with the edge $e$ in $G$. Then, if $Q$ is closed and $d_{\mathbb F}(u)=f(u)-1$, $F_2$ need not necessarily be an $f$-matching. Although $F_2$ is ostensibly not used afterwards, the condition that $F_2$ is indeed an $f$-matching is important in order to have that $\mathbf{i}_{F_1}$ is $\Pi$-extremal.
\end{remark}

\begin{remark}
In our counterexample appears a multiple edge. Note that there are also simple graphs with $\PfG\neq Q_f(G)$. For instance, consider the graph $C_4$ with a chord, where two consecutive vertices on the $4$-cycle receive the value $2$ and the others $1$.
\end{remark}

\newpage
\section{The fractional $f$-chromatic index}

Once one has Edmonds' matching polytope theorem~\ref{Edmonds1}, it is straightforward calculating to prove the following result, which was observed by Seymour \cite{Seymour} and Stahl \cite{Stahl} first.

\begin{corollary} \label{Stahl}

Any graph satisfies $$\ci^*(G)=\max\{\D^*(G),\w^*(G)\}\,\mbox{.}$$

\end{corollary}

The question is whether one can add the subscript $f$ to the three fractional graph parameters and the statement remains true. According to Zhang, Yu, and Liu \cite{Zhang}, the answer is yes. But we have already seen that their description of the $f$-matching polytope is incorrect, and we will also present a counterexample demonstrating that generally $\ciff(G)\neq\max\{\Dff(G),\wff(G)\}$. Foremost, we define the fractional graph parameter $$\Gammaff(G)=\max_{U\In V(G),\;F\In\d U,\;f(U)+|F|\geq 2}\frac{|E[U]|+|F|}{\lflr\frac{1}{2}(f(U)+|F|)\rflr}$$ (with $\Gammaff(G)=0$ if $G$ has no edge), and set $\Gammaf(G)=\lceil\Gammaff(G)\rceil$. This definition is motivated by the inequalities (c). Analogous to the derivation of Corollary~\ref{Stahl}, Theorem~\ref{Edmonds2} induces the following combinatorial characterization of the fractional $f$-chromatic index (see \cite{Stiebitz}, Theorem~B.11).

\begin{corollary} \label{Index}

Let $(G,f)$ be an arbitrary weighted graph with $\Dff(G)\geq1$. Then, $$\ciff(G)=\max\{\Dff(G),\Gammaff(G)\}\,\mbox{.}$$

\end{corollary}

In order to present the counterexample, and also after that, we will use the fact that for positive reals $a_1,\dots,a_n,b_1,\dots,b_n$, the estimate
\begin{equation} \label{max}
\frac{a_1+\cdots+a_n}{b_1+\cdots+b_n}\leq\max_{1\leq i\leq n}\frac{a_i}{b_i}
\end{equation}
holds.

\begin{example}

Consider the graph $G$ with six vertices $v_1,\cdots,v_6$. The vertices $v_1$ and $v_4$ may be connected by one edge. For a $k\in\N$, let $|E(v_1,v_2)|=|E(v_1,v_3)|=|E(v_4,v_5)|=|E(v_4,v_6)|=k$ and $|E(v_2,v_3)|=|E(v_5,v_6)|=k+1$. We therefore have $d(v)=2k+1$ for all $v\in V(G)$. The vertex function $f$ may assign the value $2$ to all vertices. Thus, $\Dff(G)=\frac{2k+1}{2}=k+\frac{1}{2}$. For every subgraph $H\In G$, $f(H)$ is even. Using (\ref{max}), we obtain $$\frac{|E(H)|}{\lfloor\frac{1}{2}f(H)\rfloor}=\frac{2|E(H)|}{f(H)}=\frac{\sum_{v\in V(H)}d_H(v)}{\sum_{v\in V(H)}f(v)}\leq\max_{v\in V(H)}\frac{d_H(v)}{f(v)}\leq\Dff(G)$$ and therefore $\wff(G)\leq\Dff(G)$. Now choose $U=\{v_1,v_2,v_3\}$ and $F=E(v_1,v_4)\In\d U$. Then we have $f(U)+|F|=3\cdot2+1=7$ and $|E[U]|+|F|=(3k+1)+1=3k+2$. Thus, $$\Gammaff(G)\geq\frac{|E[U]|+|F|}{\lflr\frac{1}{2}(f(U)+|F|)\rflr}=\frac{3k+2}{3}=k+\frac{2}{3}>k+\frac{1}{2}=\Dff(G)\,\mbox{.}$$ Corollary~\ref{Index} and $\Gammaff(G)>\Dff(G)\geq\wff(G)$ yield $$\ciff(G)=\max\{\Dff(G),\Gammaff(G)\}>\max\{\Dff(G),\wff(G)\}\,\mbox{.}$$

\end{example}

So the fractional graph parameter $\Gammaff$ seems to be more suitable for the consideration of weighted graphs than the fractional $f$-density $\wff$. We therefore tried to extend the presented counterexample to a counterexample for Conjecture~\ref{Nakano}, yet with no success. Instead, we observed that it makes no difference if one formulates this conjecture with $\wf(G)$ or $\Gammaf(G)$.

\begin{lemma} \label{mylemma}

Any weighted graph satisfies $$\max\{\Df(G)+1,\Gammaf(G)\}=\max\{\Df(G)+1,\wf(G)\}\,\mbox{.}$$

\end{lemma}

\begin{myproof} Since $\wff(G)\leq\Gammaff(G)$, we have to verify $\Gammaf(G)\leq\max\{\Df(G)+1,\wf(G)\}$ only. If $G$ has no edge, this inequality is satisfied. So assume that $E(G)\neq\emptyset$. The maximum in the definition of $\Gammaff(G)$ may be attained by $U\In V(G)$ and $F\In\d U$. With $2|E[U]|+|F|\leq2|E[U]|+|\d_GU|=\sum_{u\in U}d_G(u)$ follows $$\Gammaff(G)=\frac{|E[U]|+|F|}{\lflr\frac{1}{2}(f(U)+|F|)\rflr}\leq\frac{2(|E[U]|+|F|)}{f(U)+|F|-1}\leq\frac{\sum_{u\in U}d_G(u)+|F|}{f(U)+|F|-1}\,\mbox{.}$$ Moreover, applying (\ref{max}), we find $$\frac{\sum_{u\in U}d_G(u)}{f(U)}=\frac{\sum_{u\in U}d_G(u)}{\sum_{u\in U}f(u)}\leq\max_{u\in U}\frac{d_G(u)}{f(u)}\leq\Dff(G)\leq\Df(G)\,\mbox{.}$$

\noindent{\textbf{Case 1.} $\quad|F|\geq2$}

\noindent Then, $|F|-1$ is positive, and with (\ref{max}) and $\Df(G)\geq1$ follows $$\Gammaff(G)\leq\max\left\{\frac{\sum_{u\in U}d_G(u)}{f(U)},\frac{|F|}{|F|-1}\right\}\leq\max\{\Df(G),2\}\leq\Df(G)+1\,\mbox{.}$$
\noindent{\textbf{Case 2.} $\quad|F|=1$} $$\hspace{-35pt}\Rightarrow\qquad\Gammaff(G)\leq\frac{\sum_{u\in U}d_G(u)+1}{f(U)}\leq\Df(G)+\frac{1}{f(U)}\leq\Df(G)+1$$
\noindent{\textbf{Case 3.} $\quad|F|=0$}

\noindent Then set $H=G[U]$ and deduce $$\Gammaff(G)=\frac{|E[U]|}{\lflr\frac{1}{2}f(U)\rflr}=\frac{|E(H)|}{\lflr\frac{1}{2}f(H)\rflr}\leq\wff(G)\leq\wf(G)\,\mbox{.}$$
Consequently, $\Gammaff(G)\leq\max\{\Df(G)+1,\wf(G)\}$, and since the maximum is an integer, we have proved $\Gammaf(G)=\lceil\Gammaff(G)\rceil\leq\max\{\Df(G)+1,\wf(G)\}$.
\end{myproof}

\bigskip We can infer from Corollary~\ref{Index} that $\lceil\ciff(G)\rceil=\max\{\Df(G),\Gammaf(G)\}$, where the case $\Dff(G)<1$ needs a little extra consideration. We therefore have $$\max\{\Df(G)+1,\Gammaf(G)\}=\max\{\Df(G)+1,\lceil\ciff(G)\rceil\}\,\mbox{,}$$ since $\lceil\ciff(G)\rceil\neq\Gammaf(G)$ implies $\lceil\ciff(G)\rceil=\Df(G)\geq\Gammaf(G)$. Thus, the upper bound in Conjecture~\ref{Nakano} can be computed efficently. Should Conjecture~\ref{Nakano} additionally emerge as true, we would have $$\lceil\ciff(G)\rceil\leq\cif(G)\leq\lceil\ciff(G)\rceil+1\,\mbox{.}$$

\section*{Acknowledgements}

I am indebted to Professor Michael Stiebitz for his great support. Not only did he introduce me to the topic, but he also gave important advice and encouraged me to finish this work.

\renewcommand{\refname}{References}

\vspace{2cm}\noindent \textbf{Stefan Glock}

\medskip\noindent \emph{Institute of Mathematics\\Technische Universit\"{a}t Ilmenau\\Postfach 100565\\D-98684 Ilmenau\\Germany}

\medskip\noindent \emph{E-mail:} \verb"stefan.glock@tu-ilmenau.de"

\end{document}